\documentclass{amsproc}

\usepackage{amssymb}
\usepackage{color}

\usepackage{amsmath,amsxtra,amssymb,amsthm,amsfonts}

\newcommand{\A}{{\mathcal A}}
\newcommand{\G}{{\mathrm G}}
\newcommand{\R}{{\mathbb R}}
\newcommand{\C}{{\mathbb C}}
\newcommand{\T}{{\Bbb T}}
\newcommand{\Z}{{\Bbb Z}}
\newcommand{\M}{{\mathcal M}}
\newcommand{\St}{{\mathcal S}}
\newtheorem{theorem}{Theorem}[section]
\newtheorem{lemma}[theorem]{Lemma}

\theoremstyle{definition}

\theoremstyle{remark}
\newtheorem{remark}[theorem]{Remark}

\newcommand{\summ}{\sum\nolimits}

\allowdisplaybreaks

\def\esssup{\mathop{\mathrm{ess \, sup}}}

\numberwithin{equation}{section}
\begin{document}
\title[Harmonic analysis associated to semigroups of operators]{An invitation to harmonic analysis \\ associated with semigroups of operators}

\author{Marius Junge, Tao Mei and Javier Parcet}



\date{}

\null

\vskip-30pt

\begin{abstract}
This article is an introduction to our recent work in harmonic analysis associated with semigroups of operators, in the effort of finding a noncommutative Calder\'on-Zygmund theory for von Neumann algebras. The classical CZ theory has been traditionally developed on metric measure spaces satisfying additional regularity properties. In the lack of such metrics ---or with very little information on the metric--- Markov semigroups of operators appear to be the right substitutes of classical metric/geometric tools in harmonic analysis. Our approach is particularly useful in the noncommutative setting but it is also valid in classical/commutative frameworks.
\end{abstract}

\maketitle

\vskip-15pt

\footnote{Junge is partially supported by the NSF DMS-1201886, Mei is partially supported by the NSF DMS-0901009 and Parcet is partially supported by ERC Grant StG-256997-CZOSQP, by Spanish Grant MTM-2010-16518 and by ICMAT Severo Ochoa project SEV-2011-0087.}

\section{A sample case}

The analysis of linear operators associated to singular kernels is a central topic in harmonic analysis and partial differential equations. A large subfamily of these maps falls under the scope of Calder\'on-Zygmund theory, which exploits the relation between metric and measure in the underlying space to give sufficient conditions for $L_p$ boundedness. Namely, the H\"ormander smoothness condition for the kernel or the Calder\'on-Zygmund decomposition combine the notions of proximity in terms of the metric with that of smallness in terms of the measure. The doubling and polynomial growth conditions between metric and measure allow to extend Calder\'on-Zygmund theory to  non-Euclidean spaces. To the best of our knowledge, the existence of a metric in the underlying space is always assumed in the literature.\medskip

Let us start with a simple example which illustrates a natural framework in the lack of such metrics. Let $M_n$ be the algebra of all $n \times n$ matrices equipped with the usual trace $\mathrm{tr}$. If $A \in M_n$, consider the spectral norm $\|A\|_\infty = \sup |\langle Ae,e \rangle|$ where the supremum is taken over all unit vectors $e$ in the $n$-dimensional Hilbert space $\ell_2(n)$. In other words, the norm of $A$ as a linear operator on $\ell_2(n)$. We may also equip $M_n$ with the Schatten $p$-norm $$\|A\|_p \, = \, \big( \mathrm{tr} \, |A|^p \big)^\frac1p,$$ where $|A|=\sqrt{A^*A}$ with $A^*$ the conjugate transpose of $A$.
The Schatten $p$-classes $S_p(n)= (M_n, \|\cdot\|_p)$ $(1\leq p\leq \infty)$ behave like function $L_p$ spaces for interpolation and duality. Namely, the dual of $S_p(n)$ is isometrically  isomorphic to $S_q(n)$ when $\frac1p + \frac1q = 1$ via the duality product $\langle A,B\rangle = \mathrm{tr} (AB^*)$ and the interpolated spaces between $S_1(n)$ and $S_\infty(n)$ are isomorphic to $S_p(n)$'s with the usual interpolation parameters. All the equivalence constants are uniform in the dimension $n$.

 Let $T$ be the \emph{triangular transform}, the map sending  $A=(a_{ij})\in M_n$ to the matrix $T(A) = (\alpha_{ij}) \in M_n$ with $\alpha_{ij} = \mathrm{sgn} (i-j) a_{ij}$.  $T$ behaves as  a natural matrix analogue of the classical Hilbert transform $H$ on  the torus $\T$ or the real line $\mathbb{R}$.
  It is bounded on $S_p(n)$ uniformly in $n$ for all $1<p<\infty$ but not for $p=\infty$. Actually let $A_n=(a_{ij}) \in M_n$ with $a_{ii}=0$ and $a_{ij}=\frac1{i-j}$ for ${1\leq i\neq j\leq n}$. Kwapie\'n and Pelczy\'nski \cite{KP70} proved in 1970 that $${\|T(A_n)\|_\infty} \, \simeq \, \log n \, {\|A_n\|_{\infty}}.$$
  One may embed $M_n$ into the space of $M_n$-valued functions on the unit circle $L_\infty({\Bbb T};M_n)$ via the following embedding $$\Phi: (a_{ij}) \in M_n \mapsto \summ_{ij} a_{ij} z^{i-j} \otimes e_{ij} \in L_\infty(\T;M_n),$$ where $e_{ij}$ denote the standard matrix units. Since the map $\Phi$ is a trace preserving $*$-homomorphism, we see that $\|\Phi(A)\|_{L_p(S_p)}=\|A\|_{S_p}$ where the $L_p(S_p)$ norm is defined for matrix-valued functions on ${\T}$ as $$\|f\|_{L_p(S_p)} \, = \, \Big( \mathrm{tr} \int_{\T} |f(x)|^p dx \Big)^\frac1p.$$ Since we have $\Phi \circ T = H \circ \Phi$, it turns our that $\|T(A)\|_{S_p(n)}\leq c_p\|A\|_{S_p(n)}$ uniformly in $n$ with $c_p \simeq p^2/(p-1)$, which was already proved by Kwapie\'n/Pelczy\'nski. We refer to \cite{NR11} for more on the relation between multipliers on $S_p$ and $L_p$.

  Recall that the Hilbert transform and a large class of singular integral operators are bounded on $\mathrm{BMO}$, or at least from $L_\infty$ to $\mathrm{BMO}$. This allowed E.M. Stein and many others after him to develop and clarify Calder\'on-Zygmund theory exploiting the notion of bounded mean oscillation. The following questions for the triangular transform on matrix algebras illustrate our main goals in the context of other singular operators acting on more general von Neumann algebras

{\sc Question 1.} \emph{Can we identify a \emph{BMO}-type norm in the matrix algebras $M_n$ such that $\|T(A)\|_{\mathrm{BMO}} \leq c \, \|A\|_{\mathrm{BMO}}$  with a constant $c$ uniform in the dimension $n$}?

{\sc Question 2.} \emph{If we let $\mathrm{BMO}(n) = (M_n, \|\cdot\|_{\mathrm{BMO}})$, do we have the desired interpolation
$[\mathrm{BMO}(n), S_1(n)]_{\frac1p} \simeq \, S_p(n)$ for $1<p<\infty$ with a uniform constant}?

Let us recall the definition of the standard BMO space in the unit circle. Let us consider $f \in L_1(\T)$ and write $f_I = \frac 1{|I|} \int_I f (x) dx$ for the average of $f$ over an arc $I \subset \T$. Set $$\| f \| _{{\rm BMO}} \, = \, \sup_{I \subset \T} \frac1{|{I}|} \int_{I} |f(x) - f_{I} | dx.$$ The mean value can be replaced by the Poisson integral
\begin{eqnarray*}
\|f\| _{{\rm BMO}} & \simeq & \sup_{{I} \subset \T} \Big( \frac1{|{I}|} \int_{I} \big| f(x) - f_{I} \big|^2 dx \Big)^\frac12 \\
& \simeq & \sup_{0<t<\infty} \esssup_{x \in \T} \Big( \int_\T p_t(x-y) \Big| f(y) - \int_\T p_t(y-z) f(z) dz \Big|^2 dy \Big)^\frac12 \\ & = & \sup_{0 < t < \infty} \Big\| P_t \big| f - P_tf \big|^2 \Big\|_\infty^\frac12 \ \simeq \ \sup_{0 < t < \infty} \Big\| P_f|f|^2 - |P_t f|^2 \Big\|_\infty^\frac12
\end{eqnarray*}
with $P_t f = p_t * f$ the operator sending $\sum_k \widehat{f}(k) z^k$ to $\sum_k e^{-t|k|} \widehat{f}(k) z^k$. The above definition of BMO via the Poisson semigroup makes it possible to study BMO-type spaces in a very general framework. In fact, via the previous embedding $\Phi$ one can have a definition of BMO on $M_n$, where ``integration on intervals or balls" is not available.
For $A=(a_{ij}) \in M_n$, consider the semigroup of operators $\St = (S_t)_{t \ge 0}$ determined by $\Phi \circ S_t = P_t \circ \Phi$, which yields $${S_t}: (a_{ij}) \in M_n \mapsto ({e^{-t|i-j|}}a_{ij}) \in M_n.$$ Then, the BMO norm we propose is given by $$\|A\|_{\mathrm{BMO}_c(\St)} \, = \, \sup_{0<t<\infty} \Big\| {S_t}|A|^2 - |{S_t}A|^2 \Big\|_\infty^{\frac12}.$$ The ``c" in BMO$_c$ means column and goes back to the roots of operator space theory and noncommutative martingale  inequalities. The reason to make this distinction comes from the fact that $\|A\|_{\mathrm{BMO}_c(\St)} \nsim \|A^*\|_{\mathrm{BMO}_c(\St)}$ in general 
because of the noncommutativity of the matrix product. Later on we will see that the ``right" BMO norm for interpolation is $$\|A\|_{\mathrm{BMO}(\St)} \, = \, \max \big\{ \|A\|_{\mathrm{BMO}_c(\St)}, \|A^*\|_{\mathrm{BMO}_c(\St)} \big\}.$$

Our definition of $\mathrm{BMO}(\St)$ allows us to answer affirmatively Question 1 for the triangular transform $T$. In fact, in this particular case we obtain an isometric map since we have $$\|T (A)\|_{\mathrm{BMO}_c(\St)} \, = \, \|A\|_{\mathrm{BMO}_c(\St)}.$$ Indeed, given $A=(a_{ij}) \in M_n$ then $|A|^2=A^*A=(\sum_k \overline{a_{ki}}a_{kj})$ and
\begin{eqnarray*}
\lefteqn{\hskip-30pt \big\langle e_i, \big( S_t|A|^2-|S_tA|^2 \big) e_j \big\rangle} \\ & = & \summ_k e^{-t|i-j|} \overline{a_{ki}} a_{kj} - \summ_k e^{-t|k-j|}e^{-t|i-k|} \overline{a_{ki}} a_{kj}, \\ \lefteqn{\hskip-30pt \big\langle e_i, \big( S_t|T(A)|^2 - |S_tT(A)|^2 \big) e_j \big\rangle} \\ & = & \summ_k e^{-t|i-j|} \, {\mathrm{sgn}(k-i) \mathrm{sgn}(k-j)} \overline{a_{ki}} a_{kj} \\ & - & \summ_k e^{-t|k-j|} \, e^{-t|i-k|} \, \mathrm{sgn}(k-i) \, \mathrm{sgn}(k-j) \, \overline{a_{ki}} a_{kj}.
\end{eqnarray*}
Note that $\mathrm{sgn}(k-i) \mathrm{sgn}(k-j) \neq 1$ iff $e^{-t|i-k|}e^{-t|k-j|}=e^{-t|i-j|}$, so that $$S_t|A|^2 - |S_tA|^2 \, = \, S_t|T(A)|^2 - |S_tT(A)|^2.$$ As the proof above is purely algebraic, it is conceivable that it could work in other settings with similar algebraic structures.
For instance, one would expect a positive answer to Question 1 whenever a similar `transference' technique works, but we certainly want a general theory beyond transference. What would be a general setting to have positive answers for Questions 1 and 2?

This article is a brief (and by no means complete) introduction to the recent work on BMO spaces associated with semigroups of operators, related interpolation results and its implications in harmonic analysis. Although part of the motivation comes from noncommutative questions, our work aims to re-build harmonic analysis in the language of semigroups of operators and with assumptions merely depending on them. This follows the tradition of E.M Stein's work in Littlewood-Paley theory \cite{St} and D. Bakry's work in Riesz transforms and hypercontractivity \cite{Ba0,Ba96,BE85}. The substitution of mean value operators by semigroups of operators would make it certainly difficult to apply `local real analysis' techniques. We understand this as a challenge required to develop harmonic
analysis in a noncommutative setting.



\section{What is a framework to setup the theory?}

 From the view point of functional analysis, $L_\infty(\T)$ and the matrix algebra $M_n$  are both {\it von Neuman algebras}, in other words, weak-$*$ closed subalgebras of the space of bounded linear operators on an Hilbert space $\mathcal{B(H)}$. Both $M_n$ and $L_\infty(\T)$ come equipped with a positive linear functional $\tau$, the usual trace $\mathrm{tr}$ for $M_n$ and the integral $\int_\T \cdot \, dt$ for $L_\infty(\T)$. In both cases, $\tau$ satisfies a few natural properties. Namely, for $f,g$ living in one of these algebras we find that $\tau$ is
\begin{itemize}
\item[i)] Tracial: $\tau(fg)=\tau(gf)$,

\item[ii)] Faithful: if $f \ge 0$ and $\tau (f)=0$ then $f=0$,

\item[iii)] Lower semi-continuous: $\tau (\sup f_i)=\sup\tau (f_i)$ when $f_i\geq0$ is increasing,

\item[iv)] Semifinite: for any $f\geq 0$, there exists $0 \le g \leq f$ such that $\tau (g) < \infty$.
\end{itemize}
This led to defining {\it semifinite von Neumann algebras} as those equipped with a trace $\tau$, which is a linear (unbounded) functional satisfying (i)-(iv). Given such a pair $(\M,\tau)$ (also known as noncommutative measure space) the {\it noncommutative $L_p$ spaces} associated to it are the completion of $f\in{\mathcal M}$ with finite norm $$\|f\|_p \, = \, \big( \tau \, |f|^p \big)^\frac1p \quad \mbox{for} \quad 1 \le p < \infty,$$ where $|f|^p = (f^*f)^{p/2}$ is constructed via functional calculus. Of course, we set $L_\infty(\M)=\M$. The Schatten $p$-classes $S_p$ and the $L_p$ spaces on a semifinite measure space $(\Omega,\mu)$  are examples of noncommutative $L_p$ spaces associated with $\M = M_n$ and ${\M}=L_\infty(\Omega,\mu)$ respectively. Another basic example arises from group von Neumann algebras, which will be considered below. Every commutative semifinite von Neumann algebra is equivalent to the space $L_\infty(\Omega,\mu)$ of essentially bounded functions on some semifinite measure space $(\Omega,\mu)$. We refer to the survey paper \cite{PX03} for more information on noncommutative $L_p$ spaces.

In this general setting, there is very little information on the ``underlying space $(\Omega,\mu)$".
This is a crucial reason for the authors to pursue analogues of classical harmonic analysis without direct assumptions on the ``underlying space". As we shall explain at the end of this paper, we also expect such efforts will help in finding the geometry on the noncommutative ``$(\Omega,\mu)$".

\subsection{Markov semigroups of operators}
Semigroups of operators acting on von Neumann algebras have been widely studied in the literature, see \cite{JLX06} and the references therein. Their importance has been impressively demonstrated by the recent work of Ozawa/Popa \cite{OP10}. It also played a role in the work on Betti numbers for von Neumann algebras \cite{CS05,Pop06}. We say a semigroup of operators $\St = (S_t)_{t\geq 0}$ on a semifinite von Neumann algebra ${\mathcal M}$ is a {\it Markov semigroup} when it satisfies
\begin{itemize}
\item Unitality: $S_t(\mathbf{1}) = \mathbf{1}$ for the unit $\mathbf{1}$ of $\M$,

\item Positivity: $S_t$ are (completely) positive maps on ${\mathcal M}$,

\item Symmetry: $\tau((S_tf) g)=\tau( f (S_tg))$ for $f,g\in L_1({\mathcal M})\cap L_\infty({\mathcal M})$,

\item Weak-$*$ limits: $\tau (S_t(f)g)\rightarrow \tau (fg)$ as $t\rightarrow 0$ for $f\in {\mathcal M}, g\in L_1(\M)$.
\end{itemize}
The first and second assumptions imply that $S_t$ is a (complete) contraction in $\M$. Moreover, the first and third assumptions imply that $S_t$ is trace preserving, so that it can be extended to a contraction in $L_p(\M)$ for all $1 \le p \le \infty$. As usual a Markov semigroup always admits an infinitesimal generator $$L \, = \, - \frac{d}{d t} {S_t}_{|_{t=0}}$$ which is a closed, densely defined operator on $L_2({\mathcal M})$. We then write $ S_t=e^{-tL}$.



Markov semigroups are closely connected to Markov processes. 
In the general (noncommutative) setting, the notion of diffusion process is not yet well defined. It is however a reasonable assumption that a Markov semigroup $\St = (S_t)_{t \ge 0}$ admits a \emph{Markov dilation},
which means that there exits a larger semifinite von Neumann algebra ${\mathcal N}$, $*$-homomorphisms ${ \pi_t}:{\mathcal M}\rightarrow {\mathcal N}$ and conditional expectations ${ E_s}$ from ${\mathcal N}$ onto $\mathcal{N}_s = \bigvee_{0\leq v\leq s}\pi_v(\M)$,
such that $$\begin{array}{rcccl} & & S_{t-s} & & \\ & {\mathcal M} & \longrightarrow & {\mathcal M} \\ \pi _t & {\Big \downarrow} &  & {\Big \downarrow} & \pi _s \\ & {\mathcal N} & \longrightarrow & \mathcal{N}_s & \\ & & E_s & &\end{array}$$ is a commutative diagram $(E_{s}\pi_{t}f=\pi_{s}S_{t-s}f)$ for all $t \geq s > 0$. A Markov dilation is called {\it a.u.continuous}  if for every $2<p<\infty$, there exists a weakly dense subset $B_p$ of $L_p({\mathcal M})$ such that for any $\varepsilon>0$, there exists a projection $P_\varepsilon$ with $\tau_{\mathcal N} (1-P_\varepsilon) <\varepsilon$ and the map $t\mapsto\pi_t(S_{T-t}f)P_\varepsilon,$ is  continuous from $[0,T]$ to ${\mathcal N}$ for every $f\in B_p$, $0<T<\infty$. A Markov semigroup of operators on commutative von Neumann algebras always admits a Markov dilation \cite{Rota}. It is a recent result that a Markov semigroup of operators on finite von Neumann algebras always admits a Markov dilation \cite{JRS}. This will be relevant for interpolation purposes below. Here are some concrete examples of Markov semigroups which admit a Markov dilation

{\bf i)} Let $(M,dx)$ denote a complete Riemannian manifold with Ricci curvature bounded below by some constant $c > -\infty$. Then, the heat semigroup $\St = (S_t)_{t \ge 0}$ generated by the Laplace-Beltrami operator $\Delta$ admits a Markov dilation. We may also consider a weighted measure $d\mu(x) = \phi(x)dx$ on $M$ instead of the Riemannian measure $dx$ and the semigroup  generated by $L =\Delta + \frac {\nabla\phi}{\phi}\cdot\nabla$.

{\bf ii)} If we write (as above) $e_{ij}$ for the standard matrix units of $M_n$, the semigroup $S_t: e_{ij} \mapsto e^{-t|i-j|} e_{ij}$ which we introduced above admits a Markov dilation. More generally, consider a Hilbert space $\mathcal{H}$ and any map $b: \{1,2,\ldots,n\} \to \mathcal{H}$. Then the operators $S_t: e_{ij} \mapsto \exp(-t \|b(i)- b(j)\|_{\mathcal{H}}^2) e_{ij}$ also form a semigroup which admits a Markov dilation. This recovers the former Poisson-like semigroup by means of Schoenberg's theorem, which relates any conditionally negative length with a cocycle in the underlying group, see below. These examples are closely related to Calder\'on-Coifmann-Weiss transference method, something which we are exploring in further detail \cite{JMP3}.


{\bf iii)} Let $\mathrm{G}$ be a discrete group equipped with a conditionally negative length function $\psi: \mathrm{G} \to \mathbb{R}_+$ as defined below in Paragraph \ref{SectGvNa}. Consider the semigroup of Fourier multipliers $S_t: \lambda(g) \mapsto e^{-t \psi(g)}\lambda(g)$. Then, the existence of a dilation for $\St = (S_t)_{t \ge 0}$ follows from Ricard's paper \cite{Ri}. As an illustration, let $\mathrm{G} = \mathbb{F}_n$ be the free group with $n$ generators and consider the word length $\psi(g) = |g|$ which counts the number of letters of $g \in \mathbb{F}_n$ when it is written in reduced form. Then, the free Poisson semigroup
$S_t: \lambda(g) \mapsto e^{-t|g|} \lambda(g)$ is a Markov semigroup of operators on the group von Neumann algebra of ${\Bbb F}_n$ which admits a Markov dilation.

{\bf iv)} All the tensor products of the Markov semigroups in the examples above.

\section{BMO spaces, interpolation and duality}

Given a Markov semigroup $\St = (S_t)_{t \ge 0}$ defined on a semifinite von Neumann algebra $\M$, we now define two BMO seminorms regarding the operators $S_t$ as alternatives to the usual averaging maps. For $f \in L_2({\mathcal M})$, set
\begin{eqnarray*}
\|f\|_{{{\mathrm{BMO}_c}(\St)}} & = & \sup_{0 < t < \infty} \Big\| S_t \big| f - S_tf \big|^2 \Big\|_\infty^{\frac 12}, \\
\|f\|_{{bmo_c(\St)}} \, & = & \sup_{0 < t < \infty} \Big\| S_t|f|^2 - |S_tf|^2 \Big\|_\infty^{\frac 12}.
\end{eqnarray*}
As indicated above, the \lq\lq right" BMO seminorms arise by taking the maximum of the column-BMO norms of $f$ and $f^*$ respectively. Let us write $\mathrm{BMO}(\St)$ and $bmo(\St)$ to indicate the corresponding spaces. In order to relate both norms, we begin by introducing the gradient forms associated to $\St$. Given a semigroup of operators $\St = (S_t)_{t \ge 0} = (e^{-tL})_{t \ge 0}$ generated by $L$, define
\begin{eqnarray*}
{ \Gamma}(f,g) & = & \frac{{ L}(f^*)g+f^*{ L}(g) {- L}(f^*g)}{2}, \\ [5pt]
{ \Gamma_2}(f,g) & = & \frac{\Gamma({ L}(f),g)+\Gamma(f,{ L}(g)) -{ L}(\Gamma(f,g))}{2}.
\end{eqnarray*}
The gradient forms above are analogues of $\nabla f^* \cdot \nabla g$ and $\nabla^2f^*\cdot\nabla^2g$ respectively in terms of the generator $L$. In fact, if $L$ is the Laplace-Beltrami operator on a Riemannian manifold we obtain an identity $\Gamma(f,g)=\nabla f^* \cdot \nabla g$. For convenience it will be assumed that there is a dense $*$-algebra $\A\subset\M$ such that $L$ is well defined on $\A$ so that $\Gamma$ is well defined on $\A$. For Markov semigroups of operators, it is always true that
$\Gamma(f,f)\geq0$. We say that $\St$ satisfies the ${\Gamma_2\geq0}$ condition if $$\Gamma_2(f,f) \ge 0 \quad \mbox{for all} \quad f \in {\A}.$$ The ``$\Gamma^2\geq 0$" condition is known to hold by the heat, Ornstein-Uhlenbeck and Jacobi semigroups \cite{Ba96}, and also by all Markov semigroups of operators on group von Neumann algebras (see below for more details).  Bakry proved that $\Gamma_2\geq0 \ \Leftrightarrow\Gamma(S_tf,S_tf)\leq S_t\Gamma(f,f)\Leftrightarrow \ 2S_{t}|S_tf|^2\leq S_{2t}|f|^2+|S_{2t}f|^2$. See \cite{Ba06, CS03} for details on the gradient forms associated with semigroups of operators.

We are now in position to compare both BMO spaces introduced above via Markov semigroups of operators. Namely, if we are given a Markov semigroup $\St = (S_t)_{t \ge 0}$ which satisfies the $\Gamma_2 \ge 0$ condition, then it was proved in \cite{JM12} the following equivalence $$\|f\|_{{{\mathrm{BMO}_c}({ \St})}} \, \simeq \, \|f\|_{{bmo_c(\St)}} + \sup_{0 < t < \infty} \big\| S_tf - S_{2t}f \big\|_\infty.$$ The next result from \cite{JM12} settles the interpolation behavior of these BMO spaces.




\begin{theorem} \label{InterpTh} Let $\St = (S_t)_{t \ge 0}$ be a Markov semigroup of operators acting on a semifinite von Neumann algebra $\M$. Then, the following interpolation identities hold via the complex interpolation method for $1<p<\infty$,
\begin{itemize}
\item[i)] If $S_t$ admits a Markov dilation, then $$\big[ \mathrm{BMO}(\St),L_1({\mathcal M}) \big]_{\frac 1p} \, = \, L_p({\mathcal M}).$$

\item[ii)] Additionally, if $S_t$  admits an a.u. continuous Markov dilation then $$\hskip4pt \big[ bmo(\St),L_1({\mathcal M}) \big]_{\frac 1p}  \, = \, L_p({\mathcal M}).$$
\end{itemize}
\end{theorem}


Feffeman-Stein's $H_1$-BMO duality heavily relies on nice geometric properties of Euclidean spaces. In recent years, Auscher/McIntosh, Doung/Yan and their coauthors have obtained remarkable results in $H_1$-BMO duality theories with less geometric/metric restrictions on the underlying measure space, see \cite{ADM04, DY05} and the references therein. The probabilistic analogue of $H_1$-BMO duality has been generalized to the noncommutative case by Pisier-Xu \cite{PX97}.

We recall in Theorem \ref{DualityTh} below a recent result from \cite{M12} which settles an $H_1$-BMO duality theory that bases  merely on assumptions on the underlying semigroup of operators. This follows the line initiated by Littlewood-Paley-Stein theory. For simplicity, we shall state the duality theorem for Markov semigroups acting on commutative von Neumann algebras $\M = L_\infty(\Omega,\mu)$. Given such a semigroup $\St = (S_t)_{t \ge 0}$, recall that we replace averages by $S_t$'s and $|\nabla f|^2$ by P.A. Meyer's $\Gamma(f,f)$. Replace also $$\int_{|x-y|<t} \varphi_t(y) \, dy \frac {dt}{t^n} \, = \, \int_0^\infty \int_{B_x(t)} \varphi_t(y) dy \frac {dt}{t^n} \qquad \mbox{by} \qquad {\int_0^\infty S_t \varphi_t \, dt}.$$ Given $f \in L_1({\mathcal M})$, this motivates to introduce the square functions
\begin{eqnarray*}
{G(f)} & = & \Big( \int_0^\infty \Gamma(S_tf, S_tf) \, dt \Big)^\frac12, \\ { S(f)} & = & \Big( \int_0^\infty S_t \Gamma(S_tf, S_tf) \, dt \Big)^\frac12.
\end{eqnarray*}
Let $H_1(\St) = \big\{ f \in L_1(\M) \, : \  \|f\|_{H_1(\St)} = \|S(f)\|_{L_1({\mathcal M})} + \|f\|_{L_1({\mathcal M})} < \infty \big\}$.

\begin{theorem} \label{DualityTh}
Let $\St = (S_t)_{t \ge 0}$ be a Markov semigroup of operators acting on a semifinite von Neumann algebra $\M$ satisfying $\Gamma_2 \geq 0$. Then, the   inclusion $\mathrm{BMO}(\St)\subset (H_1(\St))^*$ is continuous. If in addition $\St$ is analytic on $L_1(\M)$, and $\|(S_t(|S_tf|^2))^\frac12\|_{L_1}\leq c\|f\|_{L_1}$ for all $f$, then $$\mathrm{BMO}(\St) = (H_1^S)^*, \quad \|S(f)\|_{L_1} \simeq \|G(f)\|_{L_1}, \quad \|\cdot\|_{bmo(\St)}\simeq \|\cdot\|_{\mathrm{BMO}(\St)}.$$
\end{theorem}
The basic examples of Markov semigroups satisfying all the assumptions in Theorem \ref{DualityTh} are heat semigroups generated by the Laplace-Beltrami operator on a complete Riemanian manifold with a nonnegative Ricci curvature, see \cite{M08} for a related result on Tent spaces and \cite{AM,CXY} for results on group von Neumann algebras and quantum tori.

\subsection{Imaginary powers}

As a first consequence of Theorem \ref{InterpTh} in harmonic analysis, the authors of \cite{JM12} also studied the $L_\infty \to \mathrm{BMO}$ boundedness of Stein's spectral multipliers $M_a$ for the subordinated Poisson semigroups $P_t = \exp(-t\sqrt L)$ which have the form $$M_a \, = \, \int_0^\infty a(t)\frac{\partial P_t}{\partial t}dt.$$ In particular, we may consider imaginary powers of the generator.

\begin{theorem} \label{Imaginary}
Let $\St = (e^{-t L})_{t \ge 0}$ be a Markov semigroup of operators acting on a semifinite von Neumann algebra $\M$ and consider the subordinated Poisson semigroup $\mathcal{P}_t = (e^{-t \sqrt{L}})_{t \ge 0}$. Then, given $u \ge 0$ the spectral multiplier $L^{iu}$ is a bounded map $$L^{iu}: L_\infty({\mathcal M}) \to \mathrm{BMO}(\mathcal{P}).$$ Therefore, if $\mathcal{P}$ admits a Markov dilation we also obtain $$\big\| L^{iu}: L_p(\M) \to L_p(\M) \big\| \, \lesssim \, \max \Big\{ p,\frac1{p-1} \Big\} u^{-|\frac12-\frac1p|} e^{|\frac {\pi u}2-\frac {\pi u}p|}.$$
\end{theorem}

Thanks to Theorem \ref{InterpTh}, the constants in Theorem \ref{Imaginary} as $u\rightarrow \infty$ are slightly    better than Stein and Cowling's, who also considered this problem \cite{Cow83,St} in the commutative case ${\M} = L_\infty(\Omega,\mu)$. For the general $M_a$'s (which correspond to Cowling's work with $\phi=\pi/2$) Theorem \ref{Imaginary} improves the $L_p$ constants in Cowling's paper from $\max\{p,\frac1{p-1}\}^{5/2}$ to the optimal ones $\max\{p,\frac1{p-1}\}$.

\section{H\"ormander-Mihlin multipliers and Riesz transforms}

Let $\G$ be a discrete group with left regular representation $\lambda: \G \to \mathcal{B}(\ell_2(\G))$ given by $\lambda(g) \delta_h = \delta_{gh}$, where the $\delta_g$'s form the unit vector basis of $\ell_2(\G)$. Write $\mathcal{L}(\G)$ for its group von Neumann algebra, the weak operator closure of the linear span of $\lambda(\G)$. Given $f \in \mathcal{L}(\G)$, consider the standard trace $\tau_\G(f) = \langle \delta_e, f \delta_e \rangle$ where $e$ denotes the identity of $\G$.  Let $L_p(\widehat{\mathbb{G}}) = L_p(\mathcal{L}(\mathrm{G}), \tau_\G)$ denote the $L_p$ space over the noncommutative measure space $(\mathcal{L}(\mathrm{G}), \tau_\G)$.  We invite the reader to check that it coincides with the usual $L_p$ over the dual group of $\G$, when $\G$ is abelian.
Any $f\in L_p(\mathcal{L}(\mathrm{G}))$ has a Fourier series $$\summ_g \widehat{f}(g) \lambda(g) \quad \mbox{with} \quad \widehat{f}(g) = \tau_\G(f \lambda(g^{-1})) \quad \mbox{so that} \quad \tau_\G(f) = \widehat{f}(e).$$
Given $m: \G \to \C$, the associated Fourier multiplier is then given by
\begin{equation} \label{FMvNa}
T_m: \summ_g \widehat{f}(g) \lambda(g) \mapsto \summ_g m(g) \widehat{f}(g) \lambda(g).
\end{equation}
Extending classical multiplier theorems to $\mathcal{L}(\G)$ was recently achieved exploiting length functions/cocycles on $\G$. This evidences some impact of cohomology theory.

\subsection{Length functions and cocycles} \label{SectGvNa}

An \emph{affine representation} or \emph{cocycle} of $\G$ is an orthogonal representation $\alpha: \G \to O(\mathcal{H})$ over a real Hilbert space $\mathcal{H}$ together with a mapping $b: \G \to \mathcal{H}$ satisfying the cocycle law $$b(gh) = \alpha_g(b(h)) + b(g).$$ In this paper we shall say that $\psi: \G \to \R_+$ is a \emph{length function} if it vanishes at the identity $e$, $\psi(g) = \psi(g^{-1})$ for all $g \in \G$ and is conditionally negative, which means that we have $$\summ_g \beta_g = 0 \Rightarrow \summ_{g,h} \overline{\beta}_g \beta_h \psi(g^{-1}h) \le 0.$$ It is easily checked that any cocycle determines a length function by the formula $\psi(g) = \|b(g)\|_\mathcal{H}^2$. According to Schoenberg's theorem \cite{Sc} the reciprocal is also true and any length function $\psi: \G \to \R_+$ determines an affine representation $(\mathcal{H}_\psi, \alpha_\psi, b_\psi)$. In fact, let $\St_\psi = (S_{\psi,t})_{t \ge 0}$ with $S_{\psi,t}: \lambda(g) \mapsto e^{-t \psi(g)} \lambda(g)$. Then Schoenberg's theorem can also be reformulated by saying that $\psi: \G \to \R_+$ is a length function iff $\St_\psi$ is Markovian, see \cite{JMP} for further details.

\subsection{H\"ormander-Mihlin multipliers}

If $\mathrm{G} = \Z^n$ the map in \eqref{FMvNa} is a Fourier multiplier on the $n$-torus. An smooth function $\widetilde{m}: \R^n \to \C$ will be called a lifting multiplier for $m$ whenever its restriction to $\Z^n$ coincides with $m$. According to de Leeuw's restriction theorem \cite{Leeuw65}, the $L_p$ boundedness of $T_m$ follows whenever there exists a lifting multiplier defining an $L_p$-bounded map in the ambient space $\R^n$. In particular, it suffices to check the H\"ormander-Mihlin smoothness condition \cite{Ho,Mi} for $1 < p < \infty$ $$\big| \partial_\xi^\beta \, \widetilde{m}(\xi) \big| \, \lesssim \, |\xi|^{-|\beta|} \quad \mbox{for all} \quad |\beta| \le \left[\frac{n}{2}\right]+1.$$ In the context of Lie groups we may find similar formulations, where the role of $\R^n$ is replaced by the corresponding Lie algebra. A fundamental goal for us is to give sufficient differentiability conditions for the $L_p$ boundedness of multipliers on the compact dual of discrete groups. Unlike for $\Z^n$, there is no standard differential structure to construct/evaluate lifting multipliers for an arbitrary discrete $\G$. The main novelty in our approach is to identify the right endpoint spaces ---intrinsic BMO's over certain semigroups--- using a broader interpretation of tangent spaces in terms of length functions and cocycles. Both H\"ormander-Mihlin and de Leeuw classical theorems are formulated in terms of the standard cocycle given by the heat semigroup. We propose the Hilbert spaces $\mathcal{H}_\psi$ as cocycle substitutes of the Lie algebra. Let us recall a fairly simple formulation ---stronger statements require more terminology--- of our cocycle form of H\"ormander-Mihlin theorem \cite{JMP}.

\begin{theorem} \label{HMTh}
Let $\G$ be a discrete group and $$T_m: \summ_g \widehat{f}(g)
\lambda(g) \mapsto \summ_g m(g) \widehat{f}(g) \lambda(g).$$ Let
$\psi$ be a length function with $\dim \mathcal{H}_\psi = n < \infty$ and $\widetilde{m}: \mathcal{H}_\psi \to \C$ such that
\begin{itemize}
\item[a)] $\widetilde{m}$ is a $\psi$-lifting of $m \!\! :$ \ $m(g) = \widetilde{m}(b_\psi(g))$,

\vskip3pt

\item[b)] $\displaystyle \big| \partial_\xi^\beta
\widetilde{m}(\xi) \big| \, \lesssim \, \min \Big\{ |\xi|^{-
|\beta|+\varepsilon}, |\xi|^{- |\beta|-\varepsilon} \Big\}$ for
$|\beta| \le \mbox{$[\frac{n}{2}]+1$}$ and some $\varepsilon > 0$.
\end{itemize}
Then, $T_m: L_p(\widehat{\mathbb{G}}) \to
L_p(\widehat{\mathbb{G}})$ is a completely bounded multiplier for
all $1 < p < \infty$.
\end{theorem}

Completely bounded means that $T_m \otimes id$ is a multiplier on $L_{p}(\widehat{\mathbb{G\times H}})$ for every discrete group $\mathrm{H}$. The additional $\varepsilon$ is a prize we pay for noncommutativity which can be removed under alternative assumptions, like
\begin{itemize}
\item[i)] $\G$ is abelian,

\item[ii)] $b_\psi(\G)$ is a lattice in $\R^n$,

\item[iii)] $\alpha_\psi(\G)$ is a finite subgroup of $O(n)$,

\item[iv)]The multiplier is $\psi$-radial, i.e. $m(g) = h(\psi(g))$.
\end{itemize}
Theorem \ref{HMTh} is a cocycle extension of the Mihlin multiplier theorem, more than merely a noncommutative form of it. Indeed, it provides new results even for finite or Euclidean groups. For instance, we may find low dimensional injective cocycles for finite groups of large cardinality, like $\mathbb{Z}_n$ or the symmetric groups $\mathrm{S}_n$, where we find injective cocycles with $\dim \mathcal{H}_\psi = 2 << n$ and $\dim \mathcal{H}_\psi = n << n!$ respectively. In the context of $\R^n$, we may use de Leeuw's compactification theorem which relates Fourier multipliers in $\R^n$ and its Bohr compactification, so that Theorem \ref{HMTh} shows that $$T_mf (x) = \int_{\R^n} \widetilde{m}(b(\xi)) \hskip1pt \widehat{f}(\xi) \hskip1pt e^{2\pi i \langle x, \xi \rangle} \, d\xi$$ is $L_p(\R^n)$-bounded for any cocycle $b: \R^n \to \R^d$ with $\widetilde{m}$ Mihlin-smooth of degree $[\frac{d}{2}]+1$. By picking suitable cocycles, this establishes a unified approach through de Leeuw's restriction/periodization theorems \cite{Leeuw65}. H\"ormander-Mihlin classical theorem corresponds to the trivial cocycle on $\R^n$. New $L_\infty \to \mathrm{BMO}$ estimates may also be given. In fact, other cocycles provide a large family of $L_p$ multipliers in $\R^n$ ---also in $\mathbb{T}^n$--- which are apparently new. As an illustration, take $$b(\xi) = (\cos 2 \pi \alpha \xi - 1, \sin 2 \pi \alpha \xi, \cos 2 \pi \beta \xi - 1, \sin 2 \pi \beta \xi)$$ for some $\alpha, \beta \in \R_+$. Theorem \ref{HMTh} shows that the restriction of a Mihlin multiplier in $\R^4$ to this \emph{donut helix} will be an $L_p$ multiplier on $\R$. It is useful to compare it with de Leeuw's periodization theorem for compactly supported multipliers. The main difference here is the irregularity obtained from choosing $\alpha/\beta$ irrational, leading to a geodesic flow with dense orbit. Hence, $m$ oscillates infinitely often with no periodic pattern. Taking for instance $\widetilde{m}(\zeta) = |\zeta|^{2\gamma}$ with $0 < \gamma < \frac12$ and smoothly truncated outside $\mathrm{B}_3(0)$, Theorem \ref{HMTh} shows that $(\sin^2(\alpha \xi) + \sin^2(\beta \xi))^\gamma$ is an $L_p$ multiplier in $\R$. These examples are certainly less standard, since the H\"ormander smoothness condition is not satisfied. With hindsight, they can be obtained via a clever combination of classical results, we invite the reader to try! However, it seems fair to say that noncommutative inspiration was required to discover such a general statement. More details can be found in \cite{JMP}.

After \cite{JMP}, length functions and cocycles have been further exploited for other families of Fourier multipliers which escape Calder\'on-Zygmund techniques on group von Neumann algebras. We refer to \cite{PR} for directional Hilbert transforms and $L_p$ convergence of Fourier series and also to \cite{JPPP,JPPPR} for hypercontractivity bounds associated to Poisson-like semigroups.



\subsection{Riesz transforms}

P.A. Meyer and D. Bakry were the first to study Riesz transforms associated with semigroups of operators. Let $\St = (S_t)_{t \ge 0}$ be a Markov semigroup of operators acting on a commutative von Neumann algebra $\M = L_\infty(\Omega,\mu)$. Assume that $\St$ is generated by $L$, which is densely defined on $\A$ and satisfies $\Gamma_2\geq 0$. Then, Bakry \cite{Ba0} proved that
\begin{equation} \label{RieszTr}
\|\Gamma^\frac12 (f,f)\|_p\simeq c_p\|L ^{\frac 12}f\|_p
\end{equation}
for all $1<p<\infty$ and any $f \in \A$. M. Junge was the first to study   the noncommutative generalization of (\ref{RieszTr}) (\cite{Ju10}).
The authors extend (\ref{RieszTr}) to group von Neumann algebras in \cite{JMP}, we omit details. For $p\geq 2$, one side of the inequality is extended in \cite{JM10} to Markov semigroups of operators on semifinite von Neumann algebras which admit a Markov dilation and satisfies the $\Gamma_2\geq0$ condition. An $L_\infty\rightarrow \mathrm{BMO}(\St)$ estimate is obtained in \cite{JM12} under the same assumptions.

On the other hand, the presence of $\varepsilon > 0$ in Theorem \ref{HMTh} excludes some central examples like the $\psi$-\emph{directional Riesz transforms} which are naturally defined for $\eta \in \mathcal{H}_\psi$ as follows (recall that $\psi(g) = \langle b_\psi(g), b_\psi(g) \rangle_\psi$) $$R_\eta \Big( \summ_g \widehat{f}(g) \lambda(g) \Big) = - i \summ_g \frac{\langle b_\psi(g), \eta \rangle_\psi}{\sqrt{\psi(g)}} \widehat{f}(g) \lambda(g).$$ Our extension of Bakry's theorem applies to $\psi$-directional Riesz transforms for any length $\psi$, even if the associated cocycle is infinite-dimensional, see \cite{JMP} for details.


\section{What are Calder\'on-Zygmund singular integrals? }

Noncommutative Calder\'on-Zygmund theory is a brand-new field, preceded by related work of Lust-Piquard and the authors  \cite{HLMP,JM10,JMP,LuP,MP,Pa0}. A major difficulty comes from the lack of appropriate analogues of certain metric/geometric properties of Euclidean spaces. Let us consider convolution kernels for the sake of clarity. We are interested in a noncommutative form of H\"ormander's condition for the kernel
$$\esssup_{x \in \R^n} \int_{|s| > 2 |x|} \big| k(s-x) - k(s) \big| \, ds \ < \ \infty.$$
The following condition is slightly weaker
$$\sup_{s > 0} \, \esssup_{x_1,x_2 \in
\mathrm{B}_s(0)} \Big| \int_{\R^n \setminus \mathrm{B}_{5 s}(0)}
\big( k(x_1-y) - k(x_2-y) \big) f(y) \, dy \Big| \ \le \ c_h \,
\|f\|_\infty,$$
but still enough to prove $L_\infty \to \mathrm{BMO}$ boundedness from $L_2$-boundedness in the usual way. Since $x_1, x_2 \notin \mathrm{supp} (f \chi_{\R^n \setminus \mathrm{B}_{5s}(0)})$, the integral above is equal to the difference $T(f \chi_{\R^n \setminus \mathrm{B}_{5s}(0)})(x_1) - T(f \chi_{\R^n \setminus \mathrm{B}_{5s}(0)})(x_2)$ for any CZ operator $T$ associated to the kernel $k$. Defining $\delta_{\R^n}(f) = f \otimes 1_{\R^n} - 1_{\R^n} \otimes f$, we may rewrite is as
\begin{equation} \label{SmoothSC}
\sup_{s > 0} \, \Big\| \big( \chi_{\mathrm{B}_s(0)} \otimes
\chi_{\mathrm{B}_s(0)} \big) \, \delta_{\R^n} \, T(f \chi_{\R^n
\setminus \mathrm{B}_{5s}(0)}) \Big\|_{L_\infty(\R^n \times \R^n)}
\ \le \ c_h \, \|f\|_\infty.
\end{equation}
This condition is much closer to its noncommutative form, but it is still unclear how to produce a family of projections in a von
Neumann algebra playing the role of the Euclidean balls centered at $0$. Our construction is based on the existence of a Markov semigroup, so that the corresponding projections will determine some sort of \lq metric\rq${}$ governing the Markov process.

Assume for simplicity that $\mathcal{S}$ is of \emph{convolution type} $S_t = (\phi_t \otimes id_\M) \circ \Delta$, where $\Delta: \M \to \M \bar\otimes \M$ is a $*$-homomorphism and $\phi_t$ is a normal state on $\M$ with density $d_t$ in $L_1^+(\M,\tau)$. If $\M$ is a Hopf algebra with comultiplication $\Delta$, each $S_t$ \lq commutes with translations\rq${}$ $\Delta \circ S_t = (S_t \otimes id_\M) \circ \Delta$. Markov semigroups of convolution type include the maps $\lambda(g) \mapsto e^{-t \psi(g)} \lambda(g)$ considered above for any length $\psi$. Of course the classical heat semigroup on $\R^n$ also falls into this category.

Let $\mathcal{S} = (S_t)_{t \ge 0}$ be a Markov semigroup of convolution type $S_t = (\phi_t \otimes id_\M) \Delta$ acting on some
noncommutative measure space $(\M, \tau)$, recall that $d_t$ stands for the density of $\phi_t = \tau(d_t \cdot)$. A \emph{weighted spectral decomposition} for $\mathcal{S}$ is a family of projections $(q_{k,t})$ in $\M$ ---indexed by $(k,t) \in \mathbb{N} \times \mathbb{R}_+$--- which are increasing in $k$ for $t$ fixed, together with a family of positive numbers $\beta_{k,t} \in \mathbb{R}_+$ such that the following conditions hold
\begin{itemize}
\item[i)] $\displaystyle \sum_{k \ge 1} \beta_{k,t} \tau(q_{k,t})
\, \le \, c_s$,

\vskip3pt

\item[ii)] $\displaystyle d_t \, \le \, c_d \, \sum_{k \ge 1}
\beta_k (q_{k,t} - q_{k-1,t})$,

\item[iii)] $\displaystyle \sum_{k \ge 1} \beta_{k,t} w_{k,t}
\tau(q_{k,t} - q_{k-1,t}) \, \le \, c_w$ \, for \, $w_{k,t} =
\big( \sum_{j \le k} \sqrt{\frac{\tau(q_{j+1,t})}{\tau(q_j,t)}}
\big)^2$,
\end{itemize}
for some absolute constants $c_d, c_s, c_w$. The first two conditions are related to Blunck/Kunstmann's analysis \cite{BK} of
non-integral CZOs, while the third condition was inspired by Tolsa's RBMO space \cite{To}. To simplify matters, let us impose the \emph{doubling property} $\tau(q_{\alpha(k),t}) \, \le c_\alpha \tau(q_{k,t})$ for some absolute constant $c_\alpha$ and some strictly increasing function $\alpha: \mathbb{N} \to \mathbb{N}$. Weighted spectral decompositions allow us to introduce a \emph{metric type} BMO. Indeed, the mean over a ball in the given metric is generalized by the convolution type formula $Q_{k,t}f = \frac{1}{\tau(q_{k,t})} \, \tau(q_{k,t} \Delta f)$, and the corresponding BMO norm is $$\|f\|_{\mathrm{BMO}_c(\mathcal{Q})} \ = \ \sup_{0 < t < \infty} \, \sup_{k \ge 1} \, \Big\| \Big( Q_{k,t} |f|^2 - |Q_{k,t} f|^2\Big)^\frac12 \Big\|_\M.$$ Taking $\|f\|_{\mathrm{BMO}(\mathcal{Q})} = \max \big\{ \|f^*\|_{\mathrm{BMO}_c(\mathcal{Q})}, \|f\|_{\mathrm{BMO}_c(\mathcal{Q})} \big\}$, we have from \cite{JMP2}
\begin{lemma}
$$\|f\|_{\mathrm{BMO}(\mathcal{S})} \ \le \ 2 \sqrt{2} \, \sqrt{c_d(c_s+c_w)}^{\null} \, \|f\|_{\mathrm{BMO}(\mathcal{Q})}.$$
\end{lemma}
Combined with Theorem \ref{InterpTh}, we may interpolate $\mathrm{BMO}(\mathcal{Q})$ with the $L_p$ scale assuming that our Markov semigroup admits a Markov dilation. The advantage is that our metric-type BMOs are adapted to produce the $L_\infty \to \mathrm{BMO}$ estimates for CZ operators, under the right algebraic analog of H\"ormander smoothness. Let us consider the derivation $\delta_\M f = f \otimes \mathbf{1}_\M - \mathbf{1}_\M \otimes f$ and set $\mathcal{L}_a/\mathcal{R}_a$ for the left/right multiplication maps $f \mapsto af/fa$. The following result provides a noncommutative analog of Calder\'on-Zygmund
extrapolation principle via H\"ormander smoothness for $T$ in terms of our \lq\lq Markov metric" above.

\begin{theorem} \label{CZTh}
Assume that
\begin{itemize}
\item $T: L_2(\M,\tau) \to L_2(\M,\tau)$ is bounded by $c_{22}$,

\vskip8pt

\item $\displaystyle \big\| \mathcal{R}_{q_{k,t} \otimes q_{k,t}}
\delta_\M T \mathcal{R}_{q_{\alpha(k),t}^\perp}: \M \to \M \otimes
\M \big\| \, \le \, c_h$ \, for all $k,t$,

\vskip3pt

\item $\displaystyle \big\| \hskip1pt \mathcal{L}_{q_{k,t} \otimes
q_{k,t}} \delta_\M T \hskip1.5pt
\mathcal{L}_{q_{\alpha(k),t}^\perp} \hskip1pt : \M \to \M \otimes
\M \big\| \, \le \, c_h$ \, for all $k,t$.
\end{itemize}
Then, we find $$\|Tf\|_{\mathrm{BMO}_\mathcal{Q}} \, \le \, \big(
2 c_{22} \sqrt{c_\alpha} + c_h \big) \, \|f\|_\infty,$$
$$\|Tf\|_{\mathrm{BMO}_\mathcal{S}} \, \le \, 2 \sqrt{2} \,
\sqrt{c_d(c_s+c_w)} \, \big( 2 c_{22} \sqrt{c_\alpha} + c_h \big)
\, \|f\|_\infty.$$
\end{theorem}

If $\mathcal{S}$ admits an a. u. continuous Markov dilation, we get the $L_p$-boundedness of $T$ by Theorem \ref{InterpTh}. This statement is not completely rigorous ---meant for Kac algebras--- but the link with \eqref{SmoothSC} is clear and gives some flavor of the theory. More general statements including nonconvolution CZO's, nondoubling settings or sufficient conditions for cb-boundedness will also be considered in \cite{JMP2}, but we will omit details here. Of course, all our discussion above applies to any $\sigma$-finite measure space $(\Omega,\mu)$ were we have no information on the metric but the space comes equipped with a Markov semigroup of operators on $L_\infty(\Omega,\mu)$.

\begin{remark}
How does classical CZ theory relates to Theorem \ref{CZTh}? Note for the heat semigroup $\St =(e^{t \Delta})_{t \geq 0}$ on ${\Bbb R}^n$, one has the following decomposition
\begin{eqnarray*}
S_tf(x)&=&\frac1{(4\pi t)^{\frac{n}{2}}} \int_{\mathbb{R}^n} \exp(-\frac {|x-y|^2}{4t}) f(y)dy\\
&=&\frac1{(4\pi t)^{\frac{n}{2}}} \int_{\R^n} \int_{\frac {|x-y|^2}{4t}}^\infty\exp(-u)du f(y)dy\\
&=&\frac{1}{\Gamma(\frac n2+1)}\int_0^\infty e^{-u}u^{\frac n2}\frac1{|\mathrm{B}_{\sqrt{4ut}}(x)|}\int_{\mathrm{B}_{\sqrt{4ut}}(x)} f(y)dydu.
\end{eqnarray*}
Thus the heat semigroup is an average of mean value operators, which explains why the classical BMOs on Euclidean spaces are equivalent to the BMOs associated with heat semigroups. Using this decomposition, we may find a weighted spectral decomposition with $(q_{k,t}, \beta_{k,t}) = (\mathrm{B}_{\sqrt{4kt}}(0), e^{-k}/(4 \pi t)^{n/2})$ which implemented in Theorem \ref{CZTh} allows us to recover the classical CZ extrapolation theorem. More general results and examples will appear in \cite{JMP2}.
\end{remark}

\section{Metric on noncommutative measure spaces}

As explained above, regarding a von Neumann algebra as a noncommutative form of $L_\infty(\Omega,\mu)$ we miss the classical geometric/metric tools used in harmonic analysis. Let us give an evidence that our approach appears to be helpful in finding a metric on the noncommutative ``$(\Omega,\mu)$".

Unital $C^*$-algebras are understood as noncommutative forms of compact Hausdorff spaces. The question we want to address is which of them play the role of ``noncommutative metric spaces"? The notion of quantum metric space goes back to M. Rieffel \cite{Ri98}. Let  ${\mathcal M}$ denote a unital $C^*$-algebra and consider a symmetric seminorm $|\cdot|_{lip}$ defined on a unital and dense $*$-subalgebra ${\mathcal A}\subset {\mathcal M}$. Suppose $|\cdot|_{lip}$  satisfies a Leibniz condition and vanishes on ${\Bbb C}\mathbf{1}$. Define a metric on the state space $S({\mathcal M})$ as $$d(\varphi,\phi) \, = \, \sup \big\{ |\varphi(a) - \phi(a)| \; | \; a \in {\mathcal A}, \ |a|_{lip} \leq 1 \big\}.$$ The triple $({\mathcal M}, {\mathcal A}, |\cdot|_{lip})$ is called a \emph{quantum metric space} if $d$ coincides with the weak-$*$ topology on  $S({\mathcal M})$. Rieffel's definition is motivated from the following example. Let ${\mathcal M}={\mathcal C}(\Omega)$ be the space of all continuous functions on a compact metric space $\Omega$, ${\mathcal A}$ the subalgebra of all Lipchitz functions on $\Omega$ and $|\cdot|_{lip}$ the usual Lipchitz norm. Then $d(\varphi,\phi)$ describes the metric on $\Omega$ for points $\varphi,\phi\in \Omega$.

Let $\St = (e^{-tL})_{t \ge 0}$ be a Markov semigroup of operators on a $C^*$-algebra ${\mathcal M}$ with its generator $L$ defined on a dense $*$-subalgebra $\A$. Let $[L^{\alpha},f]$ be the commutator of $L^{\alpha}$ and $f$.  Rieffel \cite{Ri98} checked a few model cases and proved that the triple $({C}_{{red}}^*(\Z), {\Bbb C}[\Z], \| [L^{\alpha}, \cdot \, ] \|)$ is a  quantum metric space with generator $L=\Delta$ for all $0 < \alpha \leq \frac12$. Remark 5.13 of \cite{JMP} says  that this is also true for $\frac12 < \alpha \leq 1$. Applying the noncommutative version of Bakry's theorem \eqref{RieszTr} from \cite{JM12,JMP}, the authors show in \cite{JM10,JMP} the existence of a large class of quantum metrics on von Neumann algebras of finitely generated groups with rapid decay (e.g. hyperbolic groups) satisfying the following lower estimate for some length function $\psi$ and some $\alpha > 0$ $$\inf_{|g|=k} |\psi(g)| \geq c_{\alpha} k^{\alpha}.$$ Let us now explain briefly how the noncommutative version of Bakry's theorem helps. Ozawa and Rieffel proved in \cite{OR05} that $({\mathcal M},{\mathcal A}, |\cdot|_{lip})$ is a quantum metric space iff the unit ball of $({\mathcal A}/{\Bbb C},|\cdot|_{\mathcal A})$ is relatively compact in ${\mathcal M}$. A Sobolev embedding theorem originally due to Varopoulos \cite{VSC92} in conjunction with the Ozawa/Rieffel criterium imply that $\{ f \in {\mathcal A} = {\Bbb C}[G] \, : \, \|L^{\frac 12}(f)\|_p\leq1 \}$ is relatively compact in $C^*_{red}(\G)$ for $p>\frac {4s+2}\alpha$ and $\G$ a finitely generated group with rapid decay $s$. Recall hyperbolic groups have rapid decay with $s=1$. If we combine this fact with the noncommutative analogue of \eqref{RieszTr}, we get the relative compactness of $\{f\in {\mathcal A} \, : \, \max\{\|\Gamma_\psi(f,f)\|^{\frac12},\|\Gamma_\psi(f^*,f^*)\|^{\frac12}\}\leq1\}$ for the gradient forms associated to the semigroup $\St_\psi$. This yields the following theorem.

\begin{theorem}
Let us take $$|\cdot|_{lip} \, = \, \max \Big\{ \big\| \Gamma_\psi(f,f) \big\|^{\frac12}, \big\| \Gamma_\psi(f^*,f^*) \big\|^{\frac12} \Big\},$$
for $(\G,\psi)$ as above. Then, $\big( C^*_{red}[\G], {\Bbb C}[\G], |\cdot|_{lip} \big)$ is a quantum metric space.
\end{theorem}


\bibliographystyle{amsalpha}


\enlargethispage{2cm}

\hfill \noindent \textbf{Marius Junge} \\
\null \hfill Department of Mathematics
\\ \null \hfill University of Illinois at Urbana-Champaign \\
\null \hfill 1409 W. Green St. Urbana, IL 61891. USA \\
\null \hfill\texttt{junge@math.uiuc.edu}

\vskip5pt

\hfill \noindent \textbf{Tao Mei} \\
\null \hfill Department of Mathematics
\\ \null \hfill Wayne State University \\
\null \hfill 605, W. Kirby, Detroit, MI, 48202. USA \\
\null \hfill\texttt{mei@wayne.edu}

\vskip5pt

\hfill \noindent \textbf{Javier Parcet} \\
\null \hfill Instituto de Ciencias Matem{\'a}ticas \\ \null \hfill
CSIC-UAM-UC3M-UCM \\
\null \hfill C/ Nicol\'as Cabrera 13-15. 28049, Madrid. Spain \\
\null \hfill\texttt{javier.parcet@icmat.es}
\end{document}